\documentclass[reqno,a4paper, 12pt]{amsart}

\usepackage{pictexwd,dcpic,amssymb}
\usepackage{graphicx}
\usepackage{array}
\usepackage{multirow}

\newtheorem{theorem}{Theorem}[section]
\newtheorem{proposition}[theorem]{Proposition}
\newtheorem{lemma}[theorem]{Lemma}
\newtheorem{corollary}[theorem]{Corollary}

\theoremstyle{definition}

\theoremstyle{remark}
\newtheorem{remark}[theorem]{Remark}

\numberwithin{equation}{section}

\begin{document}

\title{Singular values of principal moduli}

\author{Ja Kyung Koo}
\address{Department of Mathematical Sciences, KAIST}
\curraddr{Daejeon 373-1, Korea} \email{jkkoo@math.kaist.ac.kr}
\thanks{}

\author{Dong Hwa Shin}
\address{Department of Mathematical Sciences, KAIST}
\curraddr{Daejeon 373-1, Korea} \email{shakur01@kaist.ac.kr}
\thanks{}

\keywords{Class field theory, complex multiplication, modular and
automorphic functions.
\newline The first author was partially supported by Basic Science Research
Program through the NRF of Korea funded by MEST (2010-0001654). The
second author was partially supported by TJ Park Postdoctoral
Fellowship.}

\begin{abstract}
Let $g$ be a principal modulus with rational Fourier coefficients
for a discrete subgroup of $\mathrm{SL}_2(\mathbb{R})$ between
$\Gamma(N)$ or $\Gamma_0(N)^\dag$ for a positive integer $N$. Let
$K$ be an imaginary quadratic field. We give a simple proof of the
fact that the singular value of $g$ generates the ray class field
modulo $N$ or the ring class field of the order of conductor $N$
over $K$. Furthermore, we construct primitive generators of ray
class fields of arbitrary moduli over $K$ in terms of Hasse's two
generators.
\end{abstract}

\maketitle

\section{Introduction}

Let $\Gamma$ be a discrete subgroup of $\mathrm{SL}_2(\mathbb{R})$
commensurable with $\mathrm{SL}_2(\mathbb{Z})$. This group acts on
the complex upper half-plane
$\mathbb{H}=\{\tau\in\mathbb{C};\mathrm{Im}(\tau)>0\}$ by fractional
linear transformations, and the orbit space
$X(\Gamma)=\Gamma\backslash\mathbb{H}^*$, where
$\mathbb{H}^*=\mathbb{H}\cup\mathbb{P}^1(\mathbb{Q})$, can be given
the structure of a compact Riemann surface
(\cite[$\S$1.5]{Shimura}). When $X(\Gamma)$ is of genus zero, a
generator of the field of all meromorphic functions on $X(\Gamma)$
is called a \textit{principal modulus for $\Gamma$}.
\par
For a positive integer $N$ we denote
\begin{eqnarray*}
\Gamma(N)&=&\{\gamma\in\mathrm{SL}_2(\mathbb{Z})~;~\gamma\equiv
\left(\begin{smallmatrix}1&0\\0&1\end{smallmatrix}\right)\pmod{N}\},\\
\Gamma_1(N)&=&\{\gamma\in\mathrm{SL}_2(\mathbb{Z})~;~\gamma\equiv
\left(\begin{smallmatrix}1&\ast\\0&1\end{smallmatrix}\right)\pmod{N}\},\\
\Gamma_0(N)&=&\{\gamma\in\mathrm{SL}_2(\mathbb{Z})~;~\gamma\equiv
\left(\begin{smallmatrix}\ast&\ast\\0&\ast\end{smallmatrix}\right)\pmod{N}\},\\
\Gamma_0(N)^\dag&=&\langle\Gamma_0(N),\Phi_N\rangle,\quad\textrm{where}~
\Phi_N=\left(\begin{smallmatrix}0&-1/\sqrt{N}\\\sqrt{N}&1\end{smallmatrix}\right).
\end{eqnarray*}
Let $\Gamma$ be a discrete subgroup of $\mathrm{SL}_2(\mathbb{R})$
with $\Gamma(N)\leq\Gamma\leq\Gamma_0(N)^\dag$ for which $X(\Gamma)$
is of genus zero. Let $g$ be a principal modulus for $\Gamma$ with
rational Fourier coefficients (if any). For an imaginary quadratic
field $K$ of discriminant $d_K$ we denote
\begin{equation}\label{theta}
\theta_K=\frac{d_K+\sqrt{d_K}}{2},
\end{equation}
which generates the ring of integers $\mathcal{O}_K$ of $K$ over
$\mathbb{Z}$. Cho-Koo (\cite[Corollary 5.2]{C-K}) showed that if
$\Gamma(N)\leq\Gamma\leq\Gamma_1(N)$, then $K(g(\theta_K))$ is the
ray class field modulo $N\mathcal{O}_K$. Furthermore, Choi-Koo
(\cite[Corollary 2.5]{SC-K}) and Cho-Koo (\cite[Corollary 4.4]{C-K})
proved that if $\Gamma=\Gamma_0(N)^\dag$, then $K(g(\theta_K))$ is
the ring class field of the order of conductor $N$ in $K$, which had
been essentially done by Chen-Yui (\cite[Theorem 3.7.5(2)]{C-Y}).
Note that they used the theory of Shimura's canonical models via his
reciprocity law (\cite[$\S$6.7, 6.8]{Shimura}).
\par
In this paper, we shall first give a simple proof of the result
concerning ray class fields (Theorem \ref{main1}) by using a theorem
of Franz (\cite[Satz 1]{Franz}). On the other hand, Stevenhagen
(\cite[$\S$3, 6]{Stevenhagen}) developed a quite explicit version of
Shimura's reciprocity law. This means that we don't need to follow
the methods of Choi-Koo and Cho-Koo which are difficult of access.
And, we can give an alternative proof of the result about ring class
fields (Theorem \ref{main2}).
\par
For an imaginary quadratic field $K$ and a positive integer $N$, let
$K_{(N)}$ be the ray class field modulo $N\mathcal{O}_K$. Cho-Koo
(\cite[Corollary 5.5]{C-K}) combined Hasse's two generators of
$K_{(N)}$ by using the result of Gross-Zagier (\cite{G-Z}) and
Dorman (\cite{Dorman}) so that they obtained a primitive generator
of $K_{(N)}$ over $K$. In exactly same way we shall construct
primitive generators of ray class fields of arbitrary moduli over
$K$ (Theorem \ref{main3}).

\section{Fields of modular functions}

Let $(r_1,r_2)\in\mathbb{Q}^2-\mathbb{Z}^2$. We define the
$k^\mathrm{th}$ \textit{Fricke function} ($k=1,2,3$) (with respect
to $(r_1,r_2)$) on $\mathbb{H}$ by
\begin{equation*}
f_{(r_1,r_2)}^{(k)}(\tau)=\left\{ \begin{array}{ll}\displaystyle
-2^73^5\frac{g_2(\tau)g_3(\tau)}{\Delta(\tau)}\wp_{(r_1,r_2)}(\tau)
&
\textrm{if}~k=1\vspace{0.1cm}\\
\displaystyle\frac{g_2(\tau)^2}{\Delta(\tau)}\wp_{(r_1,r_2)}(\tau)^2
&
\textrm{if}~k=2\vspace{0.1cm}\\
\displaystyle\frac{g_3(\tau)}{\Delta(\tau)}\wp_{(r_1,r_2)}(\tau)^3 &
\textrm{if}~k=3,\\
\end{array}\right.
\end{equation*}
where
\begin{eqnarray*}
g_2(\tau)&=&60{\sum_{m,n}}^\prime\frac{1}{(m\tau+n)^4},\quad
g_3(\tau)~=~140{\sum_{m,n}}^\prime\frac{1}{(m\tau+n)^6},\quad
\Delta(\tau)~=~g_2(\tau)^3-27g_3(\tau)^2,\\
\wp_{(r_1,r_2)}(\tau)&=&\frac{1}{(r_1\tau+r_2)^2}
+{\sum_{m,n}}^\prime\bigg(
\frac{1}{(r_1\tau+r_2-m\tau-n)^2}-\frac{1}{(m\tau+n)^2}\bigg)
\end{eqnarray*}
and the sums are taken over all $(m,n)\in\mathbb{Z}^2-\{(0,0)\}$.
For simplicity we often write $f_{(r_1,r_2)}(\tau)$ instead of
$f^{(1)}_{(r_1,r_2)}(\tau)$.

\begin{proposition}\label{transfFricke}
Let $(r_1,r_2)\in\mathbb{Q}^2-\mathbb{Z}^2$.
\begin{itemize}
\item[(i)] $f_{(r_1,r_2)}^{(k)}(\tau)$ depends only on
$\pm(r_1,r_2)\pmod{\mathbb{Z}^2}$.
\item[(ii)] $f_{(r_1,r_2)}(\tau)$ satisfies the
transformation formula
\begin{equation*}
f_{(r_1,r_2)}(\tau)\circ\gamma= f_{(r_1,r_2)\gamma}(\tau)
\end{equation*}
for every $\gamma\in\mathrm{SL}_2(\mathbb{Z})$.
\end{itemize}
\end{proposition}
\begin{proof}
(i) See \cite[p.8]{Lang}.\\
(ii) See \cite[p.64]{Lang}.
\end{proof}

Let
\begin{equation*}
j(\tau)=2^63^3\frac{g_2(\tau)^3}{\Delta(\tau)}\quad(\tau\in\mathbb{H})
\end{equation*}
be the \textit{elliptic modular function}, and denote
\begin{equation*}
\mathcal{F}_1=\mathbb{Q}(j(\tau))\quad\textrm{and}\quad
\mathcal{F}_N=\mathbb{Q}(j(\tau),f_{(r_1,r_2)}(\tau))_{(r_1,r_2)\in(1/N)\mathbb{Z}^2-\mathbb{Z}^2}\quad
(N\geq2).
\end{equation*}
Note that there are relations
\begin{equation}\label{relation}
f_{(r_1,r_2)}^{(2)}(\tau)=\frac{1}{2^{8}3^{4}}
\frac{f_{(r_1,r_2)}(\tau)^2}{(j(\tau)-2^63^3)}\quad\textrm{and}\quad
f_{(r_1,r_2)}^{(3)}(\tau)=-\frac{1}{2^{9}3^{6}}\frac{f_{(r_1,r_2)}(\tau)^3}{j(\tau)(j(\tau)-2^63^3)}.
\end{equation}
We use the notations
\begin{equation*}
q=e^{2\pi i\tau}\quad\textrm{and}\quad\zeta_N=e^{2\pi
i/N}\quad(N\geq1).
\end{equation*}

\begin{proposition}\label{expansion}
\begin{itemize}
\item[(i)] We have an expansion formula
\begin{equation*}
j(\tau)=q^{-1}\prod_{n=1}^\infty(1-q^n)^{-24}\bigg(1+240\sum_{n=1}^\infty
\sigma_3(n)q^n\bigg)^3,
\end{equation*}
where $\sigma_k(n)=\sum_{d>0,d|n}d^k$
\textup{(}$k\in\mathbb{Z}$\textup{)}.
\item[(ii)] Furthermore, if $(r_1,r_2)\in\mathbb{Q}^2-\mathbb{Z}^2$, then we
get
\begin{eqnarray*}
f_{(r_1,r_2)}(\tau)&=&q^{-1}\prod_{n=1}^\infty(1-q^n)^{-24}
\bigg(1+240\sum_{n=1}^\infty\sigma_3(n)q^n\bigg)
\bigg(1-504\sum_{n=1}^\infty\sigma_5(n)q^n\bigg)\\
&&\times\bigg(1+\frac{12q^{r_1}e^{2\pi ir_2}}{(1-q^{r_1}e^{2\pi
ir_2})^2}+12\sum_{m=1}^\infty\sum_{n=1}^\infty (nq^{(m+r_1)n}e^{2\pi
ir_2n}+nq^{(m-r_1)n}e^{-2\pi ir_2n}-2nq^{mn})\bigg).
\end{eqnarray*}
\end{itemize}
\end{proposition}
\begin{proof}
See \cite[Chapter 4 $\S$1, 2]{Lang}.
\end{proof}

Hence, each function in $\mathcal{F}_N$ has a Laurent expansion with
respect to $q^{1/N}$ with coefficients in $\mathbb{Q}(\zeta_N)$,
which is called the \textit{Fourier expansion}. Furthermore,
$\mathcal{F}_N$ is a Galois extension of $\mathcal{F}_1$ with
\begin{equation*}
\mathrm{Gal}(\mathcal{F}_N/\mathcal{F}_1)\simeq\mathrm{GL}_2(\mathbb{Z}/N\mathbb{Z})/\{\pm1_2\},
\end{equation*}
whose (right) action is given as follows: For an element
$\gamma\in\mathrm{GL}_2(\mathbb{Z}/N\mathbb{Z})/\{\pm1_2\}$ we
decompose it into
\begin{equation*}
\gamma=\gamma_1\cdot\gamma_2\quad\textrm{for}~
\gamma_1=\left(\begin{smallmatrix}1&0\\0&d\end{smallmatrix}\right)~\textrm{with}~
d=\det(\gamma)\in(\mathbb{Z}/N\mathbb{Z})^*~\textrm{and any}~
\gamma_2\in\mathrm{SL}_2(\mathbb{Z}).
\end{equation*}
First, $\gamma_1$ acts by the rule
\begin{equation*}
f(\tau)=\sum_{n>-\infty} c_nq^{n/N}\mapsto f(\tau)^{\gamma_1}=
\sum_{n>-\infty} c_n^{\sigma_d}q^{n/N},
\end{equation*}
where $\sigma_d$ is the automorphism of $\mathbb{Q}(\zeta_N)$
defined by $\zeta_N^{\sigma_d}=\zeta_N^d$. And, the action of
$\gamma_2$ is given by a fractional linear transformation
(\cite[Chapter 6 Theorem 3]{Lang}).
\par
For a discrete subgroup $\Gamma$ of $\mathrm{SL}_2(\mathbb{R})$
commensurable with $\mathrm{SL}_2(\mathbb{Z})$ we denote the
corresponding modular curve by $X(\Gamma)$. In particular, if
$\Gamma=\Gamma(N)$ (respectively, $\Gamma_1(N)$, $\Gamma_0(N)$,
$\Gamma_0(N)^\dag$) for a positive integer $N$, then we simply write
$X(N)$ (respectively, $X_1(N)$, $X_0(N)$, $X_0(N)^\dag$) for
$X(\Gamma)$. Furthermore, we let $\mathbb{C}(X(\Gamma))$ be the
field of all meromorphic functions on $X(\Gamma)$, and
$\mathbb{Q}(X(\Gamma))$ be the subfield of $\mathbb{C}(X(\Gamma))$
consisting of functions with rational Fourier coefficients.

\begin{proposition}\label{Pfunction}
Let $N$ be a positive integer.
\begin{itemize}
\item[(i)] $\mathbb{C}(X(N))=\mathbb{C}\mathcal{F}_N$.
\item[(ii)] $j(N\tau)\in\mathbb{Q}(X_0(N))$.
\item[(iii)] If $N\geq2$, then $f_{(1/N,0)}(\tau)\in\mathbb{Q}(X(N))$.
\end{itemize}
\end{proposition}
\begin{proof}
(i) See \cite[Chapter 6 Theorems 1 and 2]{Lang}.\\
(ii) See \cite[Chapter 6 Theorem 5]{Lang}.\\
(iii) See  \cite[Chapter 6 Corollary 1]{Lang}.
\end{proof}

\begin{lemma}\label{Lfunction}
Let $N$ be a positive integer.
\begin{itemize}
\item[(i)] $j(\tau)j(N\tau)$, $j(\tau)+j(N\tau)\in\mathbb{Q}(X_0(N)^\dag)$.
\item[(ii)] If $N\geq2$, then $f_{(1/N,0)}^{(k)}(N\tau)\in
\mathbb{Q}(X_1(N))$ \textup{(}$k=1,2,3$\textup{)}.
\end{itemize}
\end{lemma}
\begin{proof}
(i) Observe that
\begin{equation*}
j(\tau)\circ\Phi_N=
j(\tau)\circ\left(\begin{smallmatrix}0&-1\sqrt{N}\\\sqrt{N}&0\end{smallmatrix}\right)=
j(\tau)\circ\left(\begin{smallmatrix}0&-1\\N&0\end{smallmatrix}\right)
=j(\tau)\circ\left(\begin{smallmatrix}0&-1\\1&0\end{smallmatrix}\right)\circ
\left(\begin{smallmatrix}N&0\\0&1\end{smallmatrix}\right)=j(N\tau),
\end{equation*}
and
$\Phi_N^2=\left(\begin{smallmatrix}-1&0\\0&-1\end{smallmatrix}\right)$,
from which implies that $j(\tau)j(N\tau)$ and $j(\tau)+j(N\tau)$ are
invariant via $\Phi_N$. Hence $j(\tau)j(N\tau)$ and
$j(\tau)+j(N\tau)$ belong to $\mathbb{Q}(X_0(N)^\dag)$ by
Proposition \ref{Pfunction}(ii).\\
(ii) For $\left(\begin{smallmatrix}a&b
\\c&d\end{smallmatrix}\right)\in\Gamma_1(N)$ we find that
\begin{eqnarray*}
f_{(1/N,0)}(N\tau)\circ
\left(\begin{smallmatrix}a&b\\c&d\end{smallmatrix}\right)
&=&f_{(1/N,0)}((Na\tau+Nb)/(c\tau+d))\\
&=&(f_{(1/N,0)}(\tau)\circ\left(\begin{smallmatrix}a&Nb\\c/N&d\end{smallmatrix}\right))(N\tau)\\
&=&f_{(a/N,b)}(N\tau)\quad\textrm{by Proposition \ref{transfFricke}(ii)}\\
&=&f_{(1/N,0)}(N\tau)\quad\textrm{by Proposition
\ref{transfFricke}(i)}.
\end{eqnarray*}
Hence $f_{(1/N,0)}(N\tau)$ belongs to $\mathbb{C}(X_1(N))$.
Furthermore, it has rational Fourier coefficients by Proposition
\ref{Pfunction}(iii). The same properties hold for
$f_{(1/N,0)}^{(k)}(N\tau)$ ($k=2,3$) by (\ref{relation}).
\end{proof}

\section{Shimura's reciprocity law}

For an imaginary quadratic field $K$ of discriminant of $d_K$ we let
$\theta_K$ be as in (\ref{theta}). We denote the Hilbert class field
of $K$ by $H_K$. Let $N$ be a positive integer and $\mathcal{O}$ be
the order of conductor $N$ in $K$, namely,
$\mathcal{O}=[N\theta_K,1]$. We denote by $K_{(N)}$ and
$H_\mathcal{O}$ the ray class field modulo $N\mathcal{O}_K$ and the
ring class field of $\mathcal{O}$, respectively. The following
proposition is a consequence of the theory of complex multiplication
(\cite[Chapter 10]{Lang}).

\begin{proposition}\label{CM}
Let $K$ be an imaginary quadratic field and $N$ be a positive
integer.
\begin{itemize}
\item[(i)] $K_{(N)}=K(h(\theta_K)~;~h\in\mathcal{F}_N~\textrm{is defined at
$\theta_K$})$.
\item[(ii)] If $\mathcal{O}$ is the order of conductor $N$ in $K$, then $H_\mathcal{O}=K(j(N\theta_K))$.
\item[(iii)] If $N\geq2$, then $K_{(N)}=K(j(N\theta_K),f_{(1/N,0)}^{(k)}(N\theta_K))$, where
$k=|\mathcal{O}_K^\times|/2$.
\end{itemize}
\end{proposition}
\begin{proof}
(i) See \cite[Chapter 10 Corollary to Theorem 2]{Lang}.\\
(ii) See \cite[Chapter 10 Theorem 5]{Lang}.\\
(iii) See \cite[Satz1]{Franz}.\\
\end{proof}

Let $K$ be an imaginary quadratic field. For each positive integer
$N$ we define the matrix group
\begin{equation*}
W_{N,K}=\bigg\{\begin{pmatrix}t-B_K s & -C_K
s\\s&t\end{pmatrix}\in\mathrm{GL}_2(\mathbb{Z}/N\mathbb{Z})~;~t,s\in\mathbb{Z}/N\mathbb{Z}\bigg\},
\end{equation*}
where
\begin{equation*}
\min(\theta_K,\mathbb{Q})=X^2+B_K X+C_K=
X^2-d_KX+\frac{d_K^2-d_K}{4}.
\end{equation*}
We have an explicit description of Shimura's reciprocity law
(\cite[Propositions 6.31 and 6.34]{Shimura}) due to Stevenhagen.

\begin{proposition}\label{Shimura}
Let $K$ be an imaginary quadratic field and $N$ be a positive
integer. Then $W_{N,K}$ gives rise to the surjection
\begin{equation}\label{map}
\begin{array}{ccl}
W_{N,K}&\longrightarrow&\mathrm{Gal}(K_{(N)}/H_K)\vspace{0.1cm}\\
\alpha&\mapsto&(h(\theta_K)\mapsto h^\alpha(\theta_K)~;~
\textrm{$h(\tau)\in\mathcal{F}_N$ is defined at $\theta_K$}),
\end{array}
\end{equation}
whose kernel is
\begin{equation*}
\left\{\begin{array}{ll}
\bigg\{\pm\begin{pmatrix}1&0\\0&1\end{pmatrix},~
\pm\begin{pmatrix}-2&-5\\1&2\end{pmatrix}
\bigg\} & \textrm{if}~K=\mathbb{Q}(\sqrt{-1})\vspace{0.1cm}\\
\bigg\{\pm\begin{pmatrix}1&0\\0&1\end{pmatrix},~
\pm\begin{pmatrix}-2&-3\\1&1\end{pmatrix},~
\pm\begin{pmatrix}1&3\\-1&-2\end{pmatrix}
\bigg\} & \textrm{if}~K=\mathbb{Q}(\sqrt{-3})\vspace{0.1cm}\\
\bigg\{\pm\begin{pmatrix}1&0\\0&1\end{pmatrix}\bigg\} &
\textrm{otherwise.}
\end{array}\right.
\end{equation*}
\end{proposition}
\begin{proof}
See \cite[$\S$3]{Stevenhagen}.
\end{proof}

\begin{corollary}\label{ringGalois}
Let $K$ be an imaginary quadratic field and $\mathcal{O}$ be the
order of conductor $N$ \textup{(}$\geq1$\textup{)} in $K$. Then the
map in \textup{(\ref{map})} induces an isomorphism
\begin{equation*}
\bigg\{\begin{pmatrix} t &0\\0&t
\end{pmatrix}~;~t\in(\mathbb{Z}/N\mathbb{Z})^*\bigg\}\bigg/\bigg\{\pm
\begin{pmatrix}1&0\\0&1\end{pmatrix}\bigg\}
\stackrel{\sim}{\longrightarrow}\mathrm{Gal}(K_{(N)}/H_\mathcal{O}).
\end{equation*}
\end{corollary}
\begin{proof}
See \cite[Proposition 5.3]{K-S2}.
\end{proof}

Now we can develop an analogue of Proposition \ref{CM}(i) in the
case of ring class fields.

\begin{theorem}\label{analogue}
Let $K$ be an imaginary quadratic field and $\mathcal{O}$ be the
order of conductor $N$ \textup{(}$\geq1$\textup{)} in $K$. Then
\begin{equation}\label{ringclassfield}
H_\mathcal{O}=K(h(\theta)~;~h(\tau)\in\mathbb{Q}(X_0(N))~\textrm{is
defined at $\theta_K$}).
\end{equation}
\end{theorem}
\begin{proof}
Put $R$ be the field on the right hand side of
(\ref{ringclassfield}), which is contained in $K_{(N)}$ by
Proposition \ref{CM}(i). Since $j(N\tau)\in\mathbb{Q}(X_0(N))$ by
Proposition \ref{Pfunction}(ii) and $H_\mathcal{O}=K(j(N\theta_K))$
by Proposition \ref{CM}(ii), we have the inclusion
$H_\mathcal{O}\subseteq R\subseteq K_{(N)}$. Let $h(\tau)$ be an
element of $\mathbb{Q}(X_0(N))$ which is defined at $\theta_K$. Let
$\left(\begin{smallmatrix}t&0\\0&t\end{smallmatrix}\right)\in\mathrm{GL}_2(\mathbb{Z}/N\mathbb{Z})$
with $t\in(\mathbb{Z}/N\mathbb{Z})^*$, which can be viewed as an
element of $\mathrm{Gal}(K_{(N)}/H_\mathcal{O})$ by Corollary
\ref{ringGalois}. If we decompose
$\left(\begin{smallmatrix}t&0\\0&t\end{smallmatrix}\right)=
\left(\begin{smallmatrix}1&0\\0&t^2\end{smallmatrix}\right)
\left(\begin{smallmatrix}a&b\\c&d\end{smallmatrix}\right)$ for any
$\left(\begin{smallmatrix}a&b\\c&d\end{smallmatrix}\right)\in\mathrm{SL}_2(\mathbb{Z})$,
then we get $c\equiv0\pmod{N}$ and derive that
\begin{eqnarray*}
h(\theta_K)^{\left(\begin{smallmatrix}t&0\\0&t\end{smallmatrix}\right)}&=&
h(\tau)^{\left(\begin{smallmatrix}t&0\\0&t\end{smallmatrix}\right)}(\theta_K)
\quad \textrm{by Proposition \ref{Shimura}}\\
&=&h(\tau)^{\left(\begin{smallmatrix}1&0\\0&t^2\end{smallmatrix}\right)
\left(\begin{smallmatrix}a&b\\c&d\end{smallmatrix}\right)}(\theta_K)\\
&=&h(\tau)^{\left(\begin{smallmatrix}a&b\\c&d\end{smallmatrix}\right)}(\theta_K)
\quad\textrm{because $h(\tau)$ has rational Fourier coefficients}\\
&=&h(\theta_K)\quad\textrm{by the fact}~
\left(\begin{smallmatrix}a&b\\c&d\end{smallmatrix}\right)\in\Gamma_0(N).
\end{eqnarray*}
This implies that $h(\theta_K)\in H_\mathcal{O}$; and hence
$R\subseteq H_\mathcal{O}$. Therefore, $H_\mathcal{O}=R$, as
desired.
\end{proof}

\section{Singular values of principal moduli}

\begin{lemma}\label{defined}
Let $\Gamma$ be a discrete subgroup of $\mathrm{SL}_2(\mathbb{R})$
commensurable with $\mathrm{SL}_2(\mathbb{Z})$. If
$\mathbb{C}(X(\Gamma))=\mathbb{C}(S)$ for a subset $S$ of
$\mathbb{Q}(X(\Gamma))$, then $\mathbb{Q}(X(\Gamma))=\mathbb{Q}(S)$.
\end{lemma}
\begin{proof}
See \cite[Lemma 4.1]{K-S}.
\end{proof}

\begin{lemma}\label{inclusion}
Let $g(\tau)$ be a principal modulus with rational Fourier
coefficients for a discrete subgroup $\Gamma$ of
$\mathrm{SL}_2(\mathbb{R})$ commensurable with
$\mathrm{SL}_2(\mathbb{Z})$ for which $X(\Gamma)$ is of genus zero.
For a given $\tau_0\in\mathbb{H}$, assume that $g(\tau_0)$ is an
algebraic number. If $h(\tau)\in\mathbb{Q}(X(\Gamma))$ is defined at
$\tau_0$, then $h(\tau_0)\in\mathbb{Q}(g(\tau_0))$.
\end{lemma}
\begin{proof}
Since $\mathbb{Q}(X(\Gamma))=\mathbb{Q}(g(\tau))$ by Lemma
\ref{defined}, we can express $h(\tau)=A(g(\tau))/B(g(\tau))$ for
some relatively prime $A(X),B(X)\in\mathbb{Q}[X]$. Suppose that
$B(g(\tau_0))=0$, then $A(g(\tau_0))=0$. Hence
$\min(g(\tau_0),\mathbb{Q})$ divides both $A(X)$ and $B(X)$, which
contradicts that $A(X)$ and $B(X)$ are relatively prime. Therefore,
$B(g(\tau_0))\neq0$ and $h(\tau_0)\in\mathbb{Q}(g(\tau_0))$.
\end{proof}

\begin{theorem}\label{main1}
Let $g(\tau)$ be a principal modulus with rational Fourier
coefficients for a congruence subgroup $\Gamma$ with
$\Gamma(N)\leq\Gamma\leq\Gamma_1(N)$ for an integer $N$
\textup{(}$\geq2$\textup{)}. Let $K$ be an imaginary quadratic
field. If $g(\tau)$ is defined at $\theta_K$, then
$K_{(N)}=K(g(\theta_K))$.
\end{theorem}
\begin{proof}
Since $\Gamma\leq\Gamma_1(N)\leq\Gamma_0(N)$, we get the natural
inclusion $\mathbb{Q}(X(\Gamma))\supseteq\mathbb{Q}(X_1(N))
\supseteq\mathbb{Q}(X_0(N))$. We find that
\begin{eqnarray*}
K_{(N)}&=&K(j(N\theta_K),f_{(1/N,0)}^{(k)}(N\theta_K))\quad\textrm{with
$k=|\mathcal{O}_K^\times|/2$ by Proposition \ref{CM}(iii)}\\
&\subseteq&K(g(\theta_K))\quad\textrm{by Proposition
\ref{Pfunction}(ii), Lemmas \ref{Lfunction}(ii) and \ref{inclusion}}\\
&\subseteq&K_{(N)}\quad\textrm{by Proposition \ref{CM}(i)}.
\end{eqnarray*}
Therefore, $K_{(N)}=K(g(\theta_K))$.
\end{proof}

\begin{remark}
\begin{itemize}
\item[(i)]
Unlike \cite[Corollary 5.2]{C-K} we don't use Shimura's reciprocity
law for the proof of Theorem \ref{main1}.
\item[(ii)] Kim (\cite[Remark 3.0.7]{Kim}) showed that $X_1(N)$ has genus zero if and only if $N=1,\cdots,10,12$.
There is a list of principal moduli for such $\Gamma_1(N)$ with
rational Fourier coefficients in \cite[p.161]{K-S}.
\end{itemize}
\end{remark}

\begin{lemma}\label{ringclasslemma}
Let $K$ be an imaginary quadratic field and $\mathcal{O}$ be the
order of conductor $N$ \textup{(}$\geq2$\textup{)} in $K$ such that
$H_K\subsetneq H_\mathcal{O}$. Then,
$H_\mathcal{O}=K(j(\theta_K)j(N\theta_K),j(\theta_K)+j(N\theta_K))$.
\end{lemma}
\begin{proof}
Put $a=j(\theta_K)$ and $b=j(N\theta_K)$. Let $\sigma$ be an element
of $\mathrm{Gal}(H_\mathcal{O}/K)$ which fixes both $ab$ and $a+b$.
We then derive $(a-a^\sigma)(a-b^\sigma)=
a^2-(a^\sigma+b^\sigma)a+a^\sigma b^\sigma =a^2-(a+b)a+ab=0$. If
$a=b^\sigma$, then $H_K=K(a)=K(b^\sigma)=K(b)=H_\mathcal{O}$ by
Proposition \ref{CM}(ii), which contradicts the assumption
$H_K\subsetneq H_\mathcal{O}$. We get $a=a^\sigma$; and hence
$b=b^\sigma$ from $a+b=a^\sigma+b^\sigma$. Since
$H_\mathcal{O}=K(b)$, $\sigma$ must be the unit element. Therefore,
$H_\mathcal{O}=K(ab,a+b)$.
\end{proof}

\begin{theorem}\label{main2}
Let $g(\tau)$ be a principal modulus with rational Fourier
coefficients for either $\Gamma=\Gamma_0(N)$ or $\Gamma_0(N)^\dag$
for a positive integer $N$. In the case of $\Gamma=\Gamma_0(N)^\dag$
we further assume that $H_K\subsetneq H_\mathcal{O}$. Let $K$ be an
imaginary quadratic field and $\mathcal{O}$ be the order of
conductor $N$ in $K$. If $g(\tau)$ is defined at $\theta_K$, then
$H_\mathcal{O}=K(g(\theta_K))$.
\end{theorem}
\begin{proof}
We derive that
\begin{eqnarray*}
H_\mathcal{O}&=& \left\{\begin{array}{ll}
K(j(N\theta_K))\quad\textrm{by Proposition \ref{CM}(ii)},
& \textrm{if}~\Gamma=\Gamma_0(N)\vspace{0.1cm}\\
K(j(\theta_K)j(N\theta_K),j(\theta_K)+j(N\theta_K))\quad \textrm{by
Lemma \ref{ringclasslemma}}, & \textrm{if}~
\Gamma=\Gamma_0(N)^\dag~\textrm{and}~H_K\subsetneq H_\mathcal{O}
\end{array}\right.\\
&\subseteq&K(g(\theta_K))\quad\textrm{by Proposition
\ref{Pfunction}(ii), Lemmas
\ref{Lfunction}(i) and \ref{inclusion}}\\
&\subseteq&H_\mathcal{O}\quad\textrm{by Theorem \ref{analogue}}.
\end{eqnarray*}
Therefore, $H_\mathcal{O}=K(g(\theta_K))$.
\end{proof}

\begin{remark}\label{rationalcoeff}
\begin{itemize}
\item[(i)] It is well-known that $X_0(N)$ has genus zero if and only if
$N=1,\cdots,10$, $12,13,16,18,25$. Furthermore, Helling
(\cite{Helling}) showed that $\Gamma_0(N)^\dag$ has genus zero if
and only if
$N=1,\cdots,21,23,\cdots,27,29,31,32,35,36,39,41,47,49,50,59,71$. We
have explicit formulas for principal moduli with rational Fourier
coefficients in all cases when $\Gamma_0(N)$ or $\Gamma_0(N)^\dag$
has genus zero (\cite{C-N}).
\item[(ii)]
Let $\Gamma=\Gamma_1(N)$ or $\Gamma_0(N)$ or $\Gamma_0(N)^\dag$ for
a positive integer $N$ and $h(\tau)\in\mathbb{C}(X(\Gamma))$. Since
$\left(\begin{smallmatrix}1&1\\0&1\end{smallmatrix}\right)\in\Gamma$,
$h(\tau)$ has the Fourier expansion with respect to $q$
(\cite[pp.28--29]{Shimura}). Note that $e^{2\pi i\theta_K}$ is a
real number for any imaginary quadratic field $K$. Thus, if
$h(\tau)$ has rational Fourier coefficients and is defined at
$\theta_K$, then $h(\theta_K)$ is a real algebraic number. It
follows that
\begin{equation*}
[K(h(\theta_K)):K]=\frac{[K(h(\theta_K)):\mathbb{Q}(h(\theta_K))]
\cdot[\mathbb{Q}(h(\theta_K)):\mathbb{Q}]}{[K:\mathbb{Q}]}
=[\mathbb{Q}(h(\theta_K)):\mathbb{Q}],
\end{equation*}
which implies that $\min(h(\theta_K),K)$ is a polynomial with
rational coefficients.
\end{itemize}
\end{remark}

\section {Primitive generators of ray class fields}

For a nonzero integral ideal $\mathfrak{c}$ of an imaginary
quadratic field $K$ we denote the ray class field modulo
$\mathfrak{c}$ by $K_\mathfrak{c}$. As a consequence of the theory
of complex multiplication we get the following proposition.

\begin{proposition}\label{Hasse}
Let $K$ be an imaginary quadratic field and $\mathfrak{c}$ be a
nontrivial integral ideal of $K$. Take any element $z$ in
$\mathfrak{c}^{-1}-\mathcal{O}_K$ and let $(r_1,r_2)$ be the pair of
rational numbers such that $z=r_1\theta_K+r_2$. Then we have
\begin{equation*}
K_\mathfrak{c}=K(j(\theta_K),f^{(k)}_{(r_1,r_2)}(\theta_K)),
\end{equation*}
where $k=|\mathcal{O}_K^\times|/2$.
\end{proposition}
\begin{proof}
See \cite[p.135]{Lang}.
\end{proof}

\begin{lemma}\label{jintegral}
If $\tau_0\in\mathbb{H}$ is imaginary quadratic, then $j(\tau_0)$ is
an algebraic integer.
\end{lemma}
\begin{proof}
See \cite[Chapter 5 Theorem 4]{Lang}.
\end{proof}

\begin{lemma}\label{primitive}
Let $K$ be an imaginary quadratic field of discriminant $d_K$. For
any prime $p$ greater than $|d_K|$ and any algebraic integer $w$ we
have $\mathbb{Q}(j(\theta_K),w)=\mathbb{Q}(j(\theta_K)+pw)$.
\end{lemma}
\begin{proof}
See \cite[Claim 5.6]{C-K}.
\end{proof}

\begin{remark}
Since $j(\theta_K)$ is a real algebraic integer by the definition
(\ref{theta}), Proposition \ref{expansion}(i) and Lemma
\ref{jintegral}, one can see that $\min(j(\theta_K),K)$ has integer
coefficients as in Remark \ref{rationalcoeff}(ii). Gross-Zagier
(\cite{G-Z}) and Dorman (\cite{Dorman}) showed that all prime
factors of the discriminant of $\min(j(\theta_K),K)$ are less than
or equal to $|d_K|$. By using this fact and the primitive element
theorem for a separable field extension (\cite[Theorem
51.15]{Fraleigh}), Cho-Koo obtained Lemma \ref{primitive}
\end{remark}

\begin{lemma}\label{integralZ}
Let $g(\tau)\in\mathcal{F}_N$ for a positive integer $N$. If all the
Fourier coefficients of $g(\tau)\circ\gamma$ are algebraic integers
for each $\gamma\in\mathrm{SL}_2(\mathbb{Z})$, then $g(\tau)$ is
integral over $\mathbb{Z}[j(\tau)]$.
\end{lemma}
\begin{proof}
See \cite[Chapter 2 Lemma 2.1]{K-L}.
\end{proof}

\begin{lemma}\label{Frickeintegral}
Let $(r_1,r_2)\in(1/N)\mathbb{Z}^2-\mathbb{Z}^2$ for an integer $N$
\textup{(}$\geq2$\textup{)}. Then $N^2f_{(r_1,r_2)}(\tau)$ is
integral over $\mathbb{Z}[j(\tau)]$.
\end{lemma}
\begin{proof}
We may restrict $0\leq r_1,r_2<1$ by Proposition
\ref{transfFricke}(i). One can see from Proposition
\ref{expansion}(ii) that the Fourier coefficients of
\begin{equation*}
\left\{\begin{array}{ll} f_{(r_1,r_2)}(\tau) &
\textrm{if}~r_1\neq0\vspace{0.1cm}\\
(1-e^{2\pi ir_2})^2f_{(r_1,r_2)}(\tau) & \textrm{if}~r_1=0
\end{array}\right.
\end{equation*}
are algebraic integers. Hence the Fourier coefficients of
$N^2f_{(r_1,r_2)}(\tau)$ are algebraic integer by the fact
$N=\prod_{k=1}^{N-1}(1-\zeta_N^k)$.
\par
On the other hand, for any
$\gamma=\left(\begin{smallmatrix}a&b\\c&d\end{smallmatrix}\right)\in\mathrm{SL}_2(\mathbb{Z})$
we have
\begin{equation*}
N^2f_{(r_1,r_2)}(\tau)\circ\gamma=N^2f_{(r_1,r_2)\gamma}(\tau)=
N^2f_{(\langle r_1a+r_2c\rangle,\langle
r_1b+r_2d\rangle)\gamma}(\tau)
\end{equation*}
by Proposition \ref{transfFricke}, where $\langle x\rangle$ is the
fractional part of $x\in\mathbb{R}$ in $[0,1)$. Hence the Fourier
coefficients of $N^2f_{(r_1,r_2)}(\tau)\circ\gamma$ are also
algebraic integers by the first part of the proof. Therefore,
$N^2f_{(r_1,r_2)}(\tau)$ is integral over $\mathbb{Z}[j(\tau)]$ by
Lemma \ref{integralZ}.
\end{proof}

Now we are ready to construct primitive generators of arbitrary ray
class fields over imaginary quadratic fields.

\begin{theorem}\label{main3}
Let $K$ be an imaginary quadratic field of discriminant $d_K$ and
$\mathfrak{c}$ be a nontrivial integral ideal of $K$. Take any prime
$p$ greater than $|d_K|$ and any element $z$ in
$\mathfrak{c}^{-1}-\mathcal{O}_K$. Let $(r_1,r_2)$ be the pair of
rational numbers with a denominator $N$ \textup{(}that is,
$(r_1,r_2)\in(1/N)\mathbb{Z}^2$\textup{)} such that
$z=r_1\theta_K+r_2$. Then we obtain
\begin{equation*}
K_\mathfrak{c}=K(j(\theta_K)+pN^2f^{(k)}_{(r_1,r_2)}(\theta_K)),
\end{equation*}
where $k=|\mathcal{O}_K^\times|/2$.
\end{theorem}
\begin{proof}
If $K=\mathbb{Q}(\sqrt{-1})$ or $\mathbb{Q}(\sqrt{-3})$, then
$j(\theta_K)=1728$ or $0$, respectively (\cite[p.261]{Cox}). Hence
$f^{(k)}_{(r_1,r_2)}(\theta_K)$ is a primitive generator of
$K_\frak{c}$ over $K$ by Proposition \ref{Hasse}. So we assume that
$K\neq\mathbb{Q}(\sqrt{-1})$, $\mathbb{Q}(\sqrt{-3})$ (and hence
$k=1$). Since $N^2f{(r_1,r_2)}(\tau)$ is integral over
$\mathbb{Z}[j(\tau)]$ by Lemma \ref{Frickeintegral}, its singular
value $N^2f_{(r_1,r_2)}(\theta_K)$ is an algebraic integer by Lemma
\ref{jintegral}. Therefore, we achieve the assertion by Lemma
\ref{primitive}.
\end{proof}

\bibliographystyle{amsplain}

\begin{thebibliography}{10}

\bibitem {C-Y} I. Chen and N. Yui, \textit{Singular values of Thompson series},
Groups, difference sets, and the Monster (Columbus, OH, 1993),
255--326, Ohio State Univ. Math. Res. Inst. Publ. 4, de Gruyter,
Berlin, 1996.

\bibitem {C-K} B. Cho and J. K. Koo, \textit{Constructions of ray class fields over imaginary quadratic fields and applications},
Quart. J. Math. 61 (2010), 199--216.

\bibitem {SC-K} S. Choi and J. K. Koo, \textit{On some ring class
fields by Shimura's canonical models}, Bull. Korean Math. Soc. 45
(2008), no. 4, 709--715.

\bibitem {Cox} D. A. Cox, \textit{Primes of the form $x^2+ny^2$: Fermat, Class Field, and Complex Multiplication},
A Wiley-Interscience Publication, John Wiley \& Sons, Inc., New
York, 1989.

\bibitem {C-N} J. H. Conway and S. P. Norton, \textit{Monstrous
moonshine}, Bull. Lond. Math. Soc. 11 (1979), 308--339.

\bibitem {Dorman} D. R. Dorman, \textit{Singular moduli, modular polynomials, and the
index of the closure of $\mathbb{Z}[j(\tau)]$ in
$\mathbb{Q}(j(\tau))$}, Math. Ann. 283 (1989), no. 2, 177--191.

\bibitem {Fraleigh} J. B. Fraleigh, \textit{A First Course in Abstract
Algebra}, 7th edition, Addison-Wesley Publishing Co., 2002.

\bibitem {Franz} W. Franz, \textit{Die Teilwerte der Weberschen Tau-Funktion},
J. Reine Angew. Math. 173 (1935), 60--64.

\bibitem{G-Z} B. Gross and D. Zagier, \textit{On singular moduli},
J. Reine Angew. Math. 355 (1985), 191--220.

\bibitem {Helling} H. Helling,
\textit{Note \"{u}ber das Geschlecht gewisser arithmetischer
Gruppen}, Math. Ann. 205 (1973), 173--179.

\bibitem {Kim} C. H. Kim, \textit{Arithmetic of some modular
functions}, Ph. D. Thesis, KAIST, 1999.

\bibitem {K-S} J. K. Koo and D. H. Shin, \textit{On some arithmetic properties of Siegel functions},
Math. Zeit. 264 (2010), no. 1, 137--177.

\bibitem {K-S2} J. K. Koo and D. H. Shin, \textit{Function fields of certain arithmetic curves and application},
Acta Arith. 141 (2010), no. 4, 321--334.

\bibitem {K-L} D. Kubert and S. Lang, \textit{Modular Units},
Grundlehren der mathematischen Wissenschaften 244, Spinger-Verlag,
New York-Berlin, 1981.

\bibitem {Lang} S. Lang, \textit{Elliptic Functions},
With an appendix by J. Tate, 2nd edition, Grad. Texts in Math. 112,
Spinger-Verlag, New York, 1987.

\bibitem {Shimura} G. Shimura, \textit{Introduction to the Arithmetic Theory of Automorphic Functions},
Iwanami Shoten and Princeton University Press, Princeton, N. J.,
1971.

\bibitem{Stevenhagen} P. Stevenhagen, \textit{Hilbert's 12th problem, complex multiplication and Shimura
reciprocity}, Class Field Theory-Its Centenary and Prospect (Tokyo,
1998), 161--176, Adv. Stud. Pure Math., 30, Math. Soc. Japan, Tokyo,
2001.

\end{thebibliography}

\end{document}